\documentclass[12pt]{article}
\usepackage{amsmath,amsthm,verbatim,amssymb,amsfonts,amscd,diagrams, graphics}
\theoremstyle{plain}
\newtheorem{theorem}{Theorem}[section]
\newtheorem{corollary}[theorem]{Corollary}
\newtheorem{lemma}[theorem]{Lemma}
\newtheorem{proposition}[theorem]{Proposition}

\theoremstyle{definition}
\newtheorem{definition}[theorem]{Definition}

\theoremstyle{remark}
\newtheorem{remark}[theorem]{Remark}

\newarrow{ul}---->
\newarrow{Backwards}<----

\newcommand{\td}[1]{\tilde{#1}}
\newcommand{\into}{\hookrightarrow}
\newcommand{\Z}{\mathbb{Z}}

\newcommand{\C}{\mathbb{C}}
\newcommand{\R}{\mathbb{R}}

\renewcommand{\H}{\mathbb H}

\newcommand{\mc}[1]{\mathcal{#1}}

\newcommand{\Hom}{\text{Hom}}

\newcommand{\mf}{\mathfrak}

\begin{document}

\title{Superperverse intersection cohomology: stratification (in)dependence}
\author{Greg Friedman\\Yale University}
\date{July 8, 2004}
\maketitle
typeset=\today

\begin{abstract}
Within its traditional range of perversity parameters, intersection cohomology is a topological invariant of pseudomanifolds. This is no longer true once one allows \emph{superperversities}, in which case intersection cohomology may depend on the choice of the stratification by which it is defined. Topological invariance also does not hold if one allows stratifications with  codimension one strata. Nonetheless, both errant situations arise in important situations, the former in the Cappell-Shaneson superduality theorem and the latter in any discussion of pseudomanifold bordism.  We show that while full invariance of intersection cohomology under restratification does not hold in this generality, it does hold up to restratifications that fix the the top stratum.

\end{abstract}



\tableofcontents
\section{Introduction}\label{S: intro}

It was shown by Goresky and MacPherson \cite{GM2}, the inventors of intersection cohomology theory, that for perversity parameters in the traditional range ($\bar p(2)=0$), intersection cohomology is a topological invariant of pseudomanifolds. In other words,
intersection cohomology modules are defined on stratified pseudomanifolds, but they turn out to be independent of the choice of stratification. 
This is not true, however, for \emph{superperversities}, i.e. perversities such that $\bar p(2)>0$. Nonetheless, superperverse intersection cohomology arises naturally in the study of stratified spaces and embeddings, playing a key role in the Cappell-Shaneson superduality theorem \cite{CS} (see below) and its applications \cite{GBF2, LM}. We will demonstrate that, while full topological invariance does not hold for superperverse intersection cohomology modules, they are invariant  under restratifications that fix the top stratum. 

More specifically, suppose that $X$ is an $n$-dimensional topological pseudomanifold,  $\bar p$ is a set of perversity parameters, $\mc G$ is a system of coefficients defined on a dense open sent of $X$, and $\mf X$ is a stratification of $\mf X$ such that the domain of definition of $\mc G$ contains the top stratum of $\mf X$. Then the intersection cohomology modules $I^{\bar p}_{\mf X}H^*(X;\mc G)$ are defined. If $\bar p$ is a traditional perversity, then these modules do not depend on the choice of stratification, and they are denoted $I^{\bar p}H^*(X;\mc G)$ \cite{GM2}. However, if $\bar p$ is a superperversity, intersection cohomology is not a topological invariant; different choices of stratification may result in different intersection cohomology modules. On the other hand, we show that if $\mf X$ and $\bar{\mf X}$ are two stratifications such that the singular loci $\Sigma$ and $\bar \Sigma$ of the stratifications agree, then $I^{\bar p}_{\mf X}H^*(X;\mc G)\cong I^{\bar p}_{\bar{\mf X}}H^*(X;\mc G)$. This is the conclusion of our main theorem:

\begin{theorem}[Theorem \ref{T: main}]\label{T: main2}
Let $X$ be  an $n$-dimensional topological pseudomanifold with (possibly empty)  pseudoboundary. Let $\Sigma$ be the $n-1$  skeleton of some topological stratification of $X$, and let $\mc G$ be a system of local coefficients on $X-\Sigma$. Let $\bar p$ be a traditional perversity or superperversity. Then the Deligne sheaf $\mc P^*\in D^b(X)$ is   independent of choice of stratification of $X$ subject to $\Sigma$ and hence so are the intersection cohomology modules $I^{\bar p}_{\Sigma}H^*(X;\mc G)$. 
\end{theorem}

To clarify, we say that a stratification $\mf X$ of an $n$-dimensional topological pseudomanifold is subject to $\Sigma$ if its top skeleton $X^{n-1}$ is equal to $\Sigma$. We define the pseudoboundary of a pseudomanifold to be the closure of the stratum $X^{n-1}-X^{n-2}$. More details concerning these definitions are contained in Section \ref{S: def}. It is worth noting here, however, that we do allow stratifications with $X^{n-1}\neq X^{n-2}$, generalizing the types of stratifications that are usually allowed for pseudomanifolds. In fact, even for traditional perversities, intersection cohomology is not invariant if non-trivial pseudoboundaries are allowed. 

We proceed as follows: In Section \ref{S: def} we provide the necessary basic definitions regarding pseudomanifolds, perversities, and intersection cohomology. In  Section \ref{S: not inv}, we show that different stratifications of a topological pseudomanifold may yield different intersection cohomology modules if we allow superperversities or non-empty pseudoboundaries. In fact, this can be demonstrated as a consequence of the following proposition, which show that low-dimensional superperverse intersection cohomology modules are simply the cohomology modules of the complement of the singular locus.

\begin{proposition}[Proposition \ref{P: ultrap}]
Suppose that $\bar p(k)\geq k-1$ for all $k$, $1\leq k\leq m$. Then $I^{\bar p}_{\mf X}H^*(X;\mc G)\cong H^*(X-\Sigma;\mc G)$ for $*\leq m-1$.  
\end{proposition}

As a corollary, we observe that perversities that are ``too super'' do not provide much new information:
\begin{corollary}[Corollary \ref{C: ultrap1}]
Let $X$ be an $n$-dimensional topological pseudomanifold, and suppose that $\bar p(k)\geq k-1$ for all $k$, $1\leq k\leq n$. Then $I^{\bar p}_{\mf X}H^*(X;\mc G)\cong H^*(X-\Sigma;\mc G)$. 
\end{corollary}

In Section \ref{S: break}, we investigate why the Goresky-MacPherson proof \cite{GM2} of invariance for traditional perversities does not hold in our more general setting. In Section \ref{S: main}, we provide the proof of Theorem \ref{T: main2}, indicating the necessary modifications to the Borel treatment in \cite{Bo} of the Goresky-MacPherson proof. This involves a modification of the Goresky-MacPherson axioms for intersection cohomology. Following this proof, we provide an alternative axiomatic characterization of the Deligne sheaf that allows us to recognize it as a codimension $\geq c$ intersection cohomology theory, in the sense of Habegger and Saper \cite{HS91}, with certain coefficients. This permits us to provide an alternative conclusion to the proof by invoking the topological invariance of these theories.

We close the introduction by observing that in the setting of the Cappell-Shaneson superduality theorem \cite{CS}, our main theorem is an immediate corollary. This theorem reads as follows: 

\begin{theorem}[Cappell-Shaneson]
Let $Y^n$ be a stratified pseudomanifold, and let $\mc L$ and $\mc M$ be local systems over $Y-\Sigma$ with coefficients in finitely generated $R$-modules. Let $\bar p$ and  $\bar q$ be a pair of perversities, one traditional and one a superperversity, such that $\bar p(k)+\bar q(k)=k-1$. Suppose that if $y\in \Sigma$ then the stalks $\mc H^i(\mc{I}^{\bar p}\mc C^*(Y;\mc M)_y)$ are torsion modules over $R$. Then a perfect pairing $\mc L\otimes_R \mc M\to R_{Y-\Sigma}$ and an $R$-orientation of $Y$ induce a canonical isomorphism  $\mc{ I}^{\bar q}\mc C^*(Y;\mc L)\cong R\Hom(\mc I^{\bar p}\mc C^*(Y;\mc M), \mc D^*_Y)[m]$ in the derived category $D^b(Y)$.
\end{theorem}
In this statement, $\mc{I}^{\bar p}\mc C^*$ is the Deligne intersection chain sheaf (denoted in our paper by $\mc P^*$), and $\mc D^*_Y$ is the Verdier dualizing complex on $Y$ over $R$. It follows in this context that $I^{\bar q}H^*(X;\mc L)$, which is the hypercohomology of  $\mc{ I}^{\bar q}\mc C^*(Y;\mc L)$, is determined by $\mc{ I}^{\bar p}\mc C^*(Y;\mc L)$. If this latter complex of sheaves carries a traditional perversity, its isomorphism class in the derived category $D^b(Y)$ is independent of the stratification of $Y$. 

While our theorem necessitates fixing the top stratum, it applies to a much broader range of situations.

\section{Definitions}\label{S: def}

\subsection{Some conventions}

All rings $R$ are assumed to be Noetherian, commutative, and of finite cohomological dimension. Complexes of sheaves of $R$-modules over a space $X$ should be considered as living in the bounded derived category $D^b(R_X)$ of complexes of sheaves of $R$ modules. The equality sign ``$=$'' between complexes of sheaves indicates quasi-isomorphism or, equivalently, isomorphism in the derived category. 

If $C^*$ is a complex of sheaves or modules, $C^*[m]$ represents the shifted complex $(C[m])^i=C^{i+m}$. A single sheaf or module $F$ is identified with complex having $F$ in dimension $0$ and the zero sheaf or module in all other dimensions. We sometimes emphasize this by writing $F[0]$. For example, $\R[0]$ would be a complex whose only non-trivial member is $\R$ in dimension $0$, while $\R[-m]$ would have a single $\R$ in dimension $m$. This notation is slightly counter-intuitive but conforms with the standard shift notation in derived category theory.

\subsection{Spaces}

We define our spaces by a hybrid of the definitions presented in  Goresky-MacPherson \cite{GM2} and Borel \cite{Bo}:

\begin{definition} \emph{Topological stratifications}:

\begin{itemize}
\item A \emph{$0$-dimensional topological stratified Hausdorff space} is a countable collection of points with the discrete topology.

\item An \emph{$n$-dimensional topological stratification} of a paracompact Hausdorff space $X$ consists of a filtration $\mf X$ by closed subsets $$X=X^n\supset X^{n-1}\supseteq \cdots\supseteq X^0\supseteq X^{-1}=\emptyset$$ such that each point $x\in X^{n-k}-X^{n-k-1}$ possesses a \emph{distinguished neighborhood} $N$ homeomorphic to $\R^{n-k}\times cL$, where $cL$ denotes the open cone on $L$, $cL=L\times[0,1)/(y,0)\sim(z,0)$, and $L$ is a compact Hausdorff space possessing  a $k-1$ dimensional topological stratification. Moreover, the homeomorphism $\phi: \R^{n-k}\times cL\to N$ should respect the filtration, i.e. $\phi$ takes $\R^{n-k}\times cL^j$ homeomorphically onto $N\cap X^{n-k+j+1}$ and it takes $\R^{n-k}\times \{\text{cone pt.}\}$ homeomorphically onto $N \cap X^{n-k}$. 

\item The set $X^j$ is called the \emph{$j$-skeleton} of $X$. The sets $S_j=X^j-X^{j-1}$, which are either empty or $j$-dimensional manifolds, are called the \emph{strata} of $X$. 

\item The space $L$ occurring in the definition of distinguished neighborhoods is called the \emph{link} of $x$. All points in a connected component of a stratum have the same link. 
\end{itemize} 
\end{definition}

\begin{definition}\emph{Pseudomanifolds}:
\begin{itemize}
\item A paracompact Hausdorff  space $X$ with an $n$-dimensional topological stratification is an \emph{$n$-dimensional stratified topological pseudomanifold of dimension $n$} if $X^{n-1}=X^{n-2}$ and $S_n=X-X^{n-1}=X-X^{n-2}$ is dense in $X$. In this case $X^{n-1}=X^{n-2}$ is also called the \emph{singular locus} and denoted $\Sigma$. 

\item A paracompact Hausdorff space $X$ is a \emph{topological pseudomanifold of dimension $n$} if it admits the structure of an $n$-dimension stratified topological pseudomanifold for some filtration $\mf X$. 

\end{itemize}
\end{definition}

We would like to consider spaces that satisfy the density condition on their $n$-dimensional strata  but do not necessarily have  $X^{n-1}=X^{n-2}$. These arise, for example, in the study of \emph{pseudomanifolds with boundary}, which are usually defined so that if $x\in X^{n-1}$, the boundary, then the link of $x$ consists of one point. More generally, the link of a point in $X^{n-1}$ may consist of several points. In this case, we will call the closure of $X^{n-1}-X^{n-2}$  the \emph{pseudoboundary} of $X$.  
We  make the following official definitions: 

\begin{definition}
\begin{itemize}
\item A paracompact Hausdorff  space $X$ with an $n$-dimensional topological stratification is an \emph{$n$-dimensional stratified topological pseudomanifold with pseudoboundary} if $S_n=X-X^{n-1}$ is dense in $X$. 

\item We call the closure $\overline{X^{n-1}-X^{n-2}}$ the \emph{pseudoboundary}, and it may be empty. We also call $X^{n-1}$ the \emph{singular locus} and denote it $\Sigma$. 

\item A paracompact Hausdorff space $X$ is an \emph{$n$-dimensional topological pseudomanifold with pseudoboundary} if it admits the structure of an $n$-dimension stratified topological pseudomanifold with pseudoboundary for some filtration $\mf X$.

\end{itemize}
\end{definition}

Any topological pseudomanifold is potentially a stratified topological pseudomanifold with boundary, since a pseudoboundary can always arise through the choice of stratification. To clarify the notation, we make the following remarks: In keeping with standard terminology, we will use the term ``stratified topological pseudomanifold'' only for topological pseudomanifolds stratified by filtrations with empty pseudoboundary ($X^{n-1}=X^{n-2}$). We will use the term ``topological pseudomanifold'' only for those spaces that can be stratified by such filtrations. When speaking about the more general cases for which pseudoboundaries may be allowed, we use the terms ``stratified topological pseudomanifold with pseudoboundary'' or ``topological pseudomanifold with pseudoboundary''. In the former case, a filtration is given; in the latter case one exists. In this language, we allow the possibility of empty pseudoboundaries. If we wish to emphasize cases where the pseudoboundary is, is not, or might be empty, we will speak explicitly of topological pseudomanifolds with empty, non-empty, or possibly empty pseudoboundary.  The last case will be the default. 

We will be particularly interested in stratifications of pseudomanifolds with pseudoboundary for which the singular locus has been fixed in advance. The following definition provides the language to express this concept.

\begin{definition}
Suppose that $X$ is an $n$-dimensional topological pseudomanifold with (possibly empty)  pseudoboundary and that $\Sigma$ is a closed subset of $X$. If $\mf X$ is an $n$-dimensional topological stratification of $X$ such that $X^{n-1}=\Sigma$, then we say that $\mf X$ is \emph{subject to $\Sigma$}. 
\end{definition}

Finally, we provide some notation regarding stratifications that will be of frequent use. We assume that  $X$ is an $n$-dimensional stratified topological pseudomanifold with pseudoboundary. Then recall that $S_{n-k}=X^{n-k}-X^{n-k-1}$ denotes the $n-k$ stratum of $X$. We set $U_k=X-X^{n-k}$ and let $i_k:U_k\into U_{k+1}$ and $j_k:S_{n-k}\into U_{k+1}$ be the inclusions. We also note that $U_{k+1}=U_k\cup S_{n-k}$. If $\mc S^*$ is a complex of sheaves on $X$, then $\mc S^*_k$ denotes the restriction $\mc S^*|_{U_k}$.

\subsection{Intersection cohomology}
Now that we have presented the necessary definitions regarding spaces and stratifications, we can provide the construction of intersection cohomology on a stratified topological pseudomanifold with pseudoboundary. The basic construction was first given by Goresky and MacPherson in \cite{GM2}. We include the modifications necessary to consider superperversities and pseudoboundaries. 

Intersection cohomology requires the definition of a \emph{perversity} parameter $\bar p$. This is a function $\bar p: \Z^{\geq 1}\to \Z$  satisfying the condition $\bar p(k)\leq \bar p(k+1)\leq \bar p(k)+1$. We call a perversity \emph{traditional} if $\bar p(1)=\bar p(2)=0$. These were the original perversities  introduced by Goresky and MacPherson in \cite{GM1} (in fact, $\bar p(1)$ is usually not defined as it  is unnecessary for pseudomanifold stratifications without pseudoboundary, but there is no harm in setting it equal to $0$ for consistency). A \emph{superperversity} is a perversity such that $\bar p(2)>0$ (which also implies that $\bar p(1)\geq 0$). We could also define subperversities with $\bar p(2)<0$, but these have no use in sheaf theoretic intersection homology; see Remark \ref{R: subp},  below. 

Intersection cohomology also takes as an input a local coefficient system (locally constant sheaf) defined on $X-\Sigma$. It is not necessary that these coefficients be extendable to all of $X$. We allow coefficient stalks of finitely-generated $R$ modules for some fixed Noetherian commutative ring $R$ of finite cohomological dimension 

Given an $n$-dimensional topological pseudomanifold with pseudoboundary, a traditional or super-perversity $\bar p$, and a local coefficient system $\mc G$ of $R$ modules on $X-\Sigma$, the associated Deligne sheaf $\mc P^*$ is defined. It is an element of the derived category $D^b(R_X)$ of bounded differential sheaf complexes of $R$ modules on $X$. For simplicity and since they will remain fixed in any given discussion, we omit $X$, $\mf X$, $\bar p$, and $\mc G$ from the notation $\mc P^*$. Let $U_k=X-X^{n-k}$, and let $i_k: U_k\to U_{k+1}=U_k\cup S_{n-k}$ denote the inclusion. 
Then $\mc P^*$ is defined inductively as follows: On $U_1$, $\mc P^*_1=\mc G[0]$, the local coefficient sheaf treated as a complex of sheaves  with the only non-trivial member of the complex being $\mc G$ in dimension $0$. For $k\geq 1$, let $\mc P^*_{k+1}=\tau_{\leq \bar p(k)}Ri_{k*}\mc P^*_k$, where $\tau_{\leq \bar p(k)}$ is the truncation functor and $Ri_{k*}$ is the right derived functor of the pushforward $i_{k*}$. Then $\mc P^*=\mc P_{n+1}^*$. 

For a fixed $X$, $\mf X$, $\bar p$, and $\mc G$, the \emph{intersection cohomology} module $I^{\bar p}_{\mf X}H^i(X;\mc G)$ is defined to be the hypercohomology $\H^i(\mc P^*)$. 
Since we will show below that intersection cohomology depends only on $X$, $\bar p$, $\mc G$, and $\Sigma$ ( \emph{not} on the entire stratification $\mf X$), we will also use $I^{\bar p}_{\Sigma}H^*(X;\mc G)$ to denote intersection cohomology as defined with respect to a stratification $\mf X$ subject to $\Sigma$.

\begin{remark}
Intersection \emph{homology}, as it is usually defined via PL or singular chains (see \cite{GM1,Bo, Ki,GBF10}), is related to intersection cohomology by the formula $I^{\bar p}H_*(X;\mc G)\cong I^{\bar p}H^{n-*}(X; \mc G\otimes_R \mc O)$, where $\mc O$ is an orientation sheaf for $X-\Sigma$. Care should be taken, however, since there may be many choices for $\mc O$ if $X-\Sigma$ is disconnected. 
\end{remark}

\begin{remark}\label{R: subp} It follows immediately from the definitions that there can be no interest in studying \emph{subperversities}, those for which $\bar p(2)<0$, and hence also $\bar p(1) <0$. In this case, $\mc P_2^*=0$, and it follows that $\mc P^*=0$ unless $X$ is a manifold and $\mf X$ is the trivial stratification, in which case $\mc P^*=\mc G^*$. 
\end{remark}

\section{Superperverse intersection homology is not stratification invariant}\label{S: not inv}

In this section, we demonstrate that superperverse intersection cohomology and intersection cohomology allowing pseudoboundaries are not topological invariants - they depend upon the stratification. We begin with the following proposition, which shows that, in a certain range, superperverse intersection cohomology is simply the cohomology of the top stratum.

\begin{proposition}\label{P: ultrap}
Suppose that $\bar p(k)\geq k-1$ for all $k$, $1\leq k\leq m$. Then $I^{\bar p}_{\mf X}H^*(X;\mc G)\cong H^*(X-\Sigma;\mc G)$ for $*\leq m-1$. 
\end{proposition}

Before proving the proposition, we note some ramifications.

It is immediate from Proposition \ref{P: ultrap} that if $\bar p(k)\geq k-1$ for any $k\geq 1$ (and so $\bar (j)\geq j-1$ for all $j\leq k$), then $I^{\bar p}H^*(X;\mc G)$ depends on the singular set $\Sigma$. In particular, this will occur any time we allow  non-empty pseudoboundaries since in this case $\bar p(1)=0=1-1$. Even if we disallow pseudoboundaries, the same issues occur whenever $\bar p(2)\geq 1$, as we may still alter $\Sigma=X^{n-2}$ by restratification.  

We can provide some simple illustrations by considering nonstandard stratifications of manifolds. For example, the proposition tells us that if we stratify the  sphere $S^n$ by $S^n=X^n\supset X^0=x$ for some $x\in S^n$, then if $\bar p(k)=k-1$ for all $k$, $I^{\bar p}H^*(S^n;\R_{S^n})=H^*(S^n-x;\R)=H^*(\R^n;\R)=\R[0]$. However, if we stratify $S^n$ with the trivial filtration, then $I^{\bar p}H^*(S^n;\R_{S^n})=H^*(S^n;\R)=\R[0]\oplus \R[-n]$. Similarly, if we stratify $S^{n}$ by $S^n=X^n\supset X^{n-1}=S^{n-1}$ for the standard embedding $S^{n-1}\into S^n$, then even for the more traditional perversity $\bar p\equiv 0$, we have $\bar p(k)\geq k-1$ for the only  relevant value of $k$, $k=1$. Thus in this case $I^{\bar p}H^*(S^n;\R_{S^n})=H^*(S^n-S^{n-1};\R)=\R[0]\oplus \R[0]$. In this case we see dependence on the choice of pseudoboundary. 

Another consequence of Proposition \ref{P: ultrap} is the following corollary, which shows that  perversities that are too super do not yield interesting intersection cohomology modules.

\begin{corollary}\label{C: ultrap1}
Let $X$ be an $n$-dimensional topological pseudomanifold, and suppose that $\bar p(k)\geq k-1$ for all $k$, $1\leq k\leq n$. Then $I^{\bar p}_{\mf X}H^*(X;\mc G)\cong H^*(X-\Sigma;\mc G)$. 
\end{corollary}
\begin{proof}
Since $X$ has dimension $n$, only the values of $\bar p(k)$ for $k\leq n$ have any relevance for defining the intersection cohomology of $X$. Therefore, since $\bar p(n)\geq n-1$, we can extend $\bar p$ so that $\bar p(n+1)\geq n$ without affecting the intersection cohomology. The corollary now follows  from Proposition \ref{P: ultrap}.
\end{proof}

\begin{proof}[Proof of Proposition \ref{P: ultrap}]
Recall that  $I^{\bar p}_{\mf X}H^j(X;\mc G)= \H^j(X;\mc P^*)$, where $\mc P^*$ is the Deligne sheaf associated to $X$, $\mf X$, $\bar p$, and $\mc G$. We first claim that  it suffices to show that  if $i$ is the inclusion $i:X-\Sigma\into X$ and $\bar p(k)\geq k-1$ for  $1\leq k\leq m$, then   $ \mc P^*$ and $Ri_*\mc G$ are quasi-isomorphic up through dimension $m-1$. For suppose that this is so. Then it follows from the hypercohomology spectral sequences that $\H^j(X;Ri_*\mc G)\cong \H^j(X;\mc P^*)$ for $j\leq m-1$ (see \cite[Theorem IV.2.2]{Br}). But  $\H^j(X;\mc P^*)= I^{\bar p}_{\mf X}H^j(X;\mc G)$, while $\H^j(X;Ri_*\mc G)=\H^j(X-\Sigma;\mc G)=H^j(X-\Sigma;\mc G)$, since $\mc G$ is a local system on the manifold  $X-\Sigma$. 

We will first show that  $ \mc P^*_{k+1}=Ri_{k*}\cdots Ri_{1*}\mc G$ for $k\leq m$.  To do this, it suffices to demonstrate that  $ \tau_{\leq \bar p(k)}Ri_{k*}\mc P^*_k=Ri_{k*}\cdots Ri_{1*}\mc G$  for each $k\leq m$. Since the equality is in the derived category, we need only show that for each $x\in U_{k+1}$, $\mc H^*(Ri_{k*} \cdots Ri_{1*}\mc G)_x=0$ for $*>\bar p(k)\geq k-1$. We will proceed again by induction.

We begin  with $Ri_{1*}\mc G$ on $U_2$. We need to see that for all $x\in U_2$, $\mc H^*( Ri_{1*}\mc G)_x=0$ for $*>\bar p(1)\geq 0$. Now $U_2=U_1\cup S_{n-1}$, and $ (Ri_{1*}\mc G)|_{U_1}=\mc G$. So at each point $x\in U_1$,  $\mc H^*( Ri_{1*}\mc G)_x=H^*(\mc G_x)=G[0]$, where $G$ is the stalk of $\mc G$ at $x$. Next consider $ S_{n-1}$. If $S_{n-1}$ is empty, there is nothing to show. If $S_{n-1}$ is not empty and $x\in S_{n-1}$, then the link of $x$ consists of a finite number of points, $L\cong \amalg y_i$, and $x$ has a fundamental system of distinguished neighborhoods $N$ of the form $\R^{n-1}\times cL$. 
In this case, $\mc H^*( Ri_{1*}\mc G)_x\cong \lim_{x\in N}\H^*(N; Ri_{1*}\mc G)$, where the limit is taken over distinguished neighborhoods of $x$. But by Lemma V.3.9 of \cite{Bo},  $\H^*(N; Ri_{1*}\mc G)\cong \H^*(L;\mc G|_L)=H^*(L;\mc G|_L)$. This lemma is stated for stratified pseudomanifolds but the proof holds just as well for stratified pseudomanifolds with pseudoboundary. By further results in \cite[\S 3]{Bo}, it even follows that the direct system $\H^*(N; Ri_{1*}\mc G)$ over distinguished neighborhoods $N$ is essentially constant, but the resulted already quoted is sufficient to demonstrate that  $\mc H^*( Ri_{1*}\mc G)_x=0$ for $*>0$. 

Now suppose inductively that we have shown for all $j$, $1\leq j<k$, that  for $x\in U_{j+1}$, $\mc H^*(Ri_{j*} \cdots Ri_{1*}\mc G)_x=0$ for $*>j-1$. We consider $Ri_{k*} \cdots Ri_{1*}\mc G$ on  $U_{k+1}=U_k\cup S_{n-k}$. Since $ (Ri_{k*} \cdots Ri_{1*}\mc G)|_{U_k}= Ri_{k-1*} \cdots Ri_{1*}\mc G $, we already know by induction that if $x\in U_k$ then $\mc H^*(Ri_{k*} \cdots Ri_{1*}\mc G)_x=0$ for $*>k-2$. Next we consider points in $S_{n-k}$. If $S_{n-k}$ is empty, then there is nothing to show and this induction step is finished. Suppose then that $S_{n-k}\neq \emptyset$ and $x\in S_{n-k}$. Once again,  Lemma V.3.9 of \cite{Bo} tells us that  $\mc H^*(Ri_{k*}\cdots Ri_{1*}\mc G)_x\cong \H^*(L; (Ri_{k-1*}\cdots Ri_{1*}\mc G)|_L)$, where $L$ is the link of $x$. To invoke this lemma, it is only necessary to apply part (b) of the same lemma inductively to see that $Ri_{k-1*}\cdots Ri_{1*}\mc G$ is $\mf X$-cohomologically locally constant ($\mf X$-clc). But $ Ri_{k-1*}\cdots Ri_{1*}\mc G = R(i_{k-1*}\cdots Ri_{1*})\mc G $, so $\H^*(L; (Ri_{k-1*}\cdots Ri_{1*}\mc G)|_L)=\H^*(L-L^{k-2};\mc G)=H^*(L-L^{k-2};\mc G)$. Since $L-L^{k-2}$ is a $k-1$ manifold, its cohomology with coefficients in the local system $\mc G$ is $0$ for $*>k-1$.

Up to this point, we have established that $\mc P^*_{m+1}=Ri_{m*} \cdots Ri_{1*}\mc G$. So from here it suffices to show that if $\mc S^*$ is any complex of sheaves defined on $U_{m+1}$, then 
$$\mc H^j(Ri_{n*} \cdots Ri_{(m+1)*}\mc S^*) =\mc H^j(\tau_{\leq \bar p(n)}Ri_{n*} \cdots \tau_{\leq \bar p(m+1)}Ri_{(m+1)*} \mc S^* )$$ 
for $j\leq  m-1$.  But since $\bar p$ is a perversity and $\bar p(k)\geq k-1$ for $1\leq k\leq m$, we must have that $\bar p(k)\geq m-1$ for $k>m$. So it suffices to show that if there is a morphism $\phi: \mc A^*\to \mc B^*$ between two $\mf X$-clc sheaves on $U_k$ for some $k\geq m+1$ and if $\phi$ is a  quasi-isomorphism up to dimension $m-1$, then the induced map $Ri_{k*}\mc A^*\to \tau_{\leq \bar p(k)}Ri_{k*}\mc B^*$ is also a quasi-isomorphism up to dimension $m-1$. For this we need only show that $Ri_{k*}\mc A^*$ and $Ri_{k*}\mc B^*$ are quasi-isomorphic up to dimension $m-1$ since clearly $Ri_{k*}\mc B^*$ and $\tau_{\leq \bar p(k)}Ri_{k*}\mc B^*$ are quasi-isomorphic up to dimension $\bar p(k)\geq m-1$. 

To check this last condition, we note that for $x\in U_k$ and any $\mf X$-clc sheaf complex $\mc S^*$, $(Ri_{k*}\mc S^*)_x=\mc S^*_x$, so by hypothesis $Ri_{k*}\mc A^*$ and $Ri_{k*}\mc B^*$ are quasi-isomorphic up to dimension $m-1$ at each point $x\in U_k$. If $x\in S_{n-k}=U_{k+1}-U_k$, then again by \cite[Lemma V.3.9.a]{Bo}, $H^j(Ri_{k*}\mc S^*)_x=\H^j(L;\mc S^*)$. So we must show that the map $\H^j(L;\mc A^*|_L)\to \H^j(L;\mc B^*|_L)$ induced by $\phi$ is a quasi-isomorphism up to dimension $m-1$. But once more employing \cite[Theorem IV.2.2]{Br}, this is true if $\mc A^*|_L$ and $\mc B^*|_L$ are quasi-isomorphic up through dimension $m-1$, and this holds by assumption since $L\subset U_k$. 

This concludes the proof of the theorem.
\end{proof}

Another consequence of  Proposition \ref{P: ultrap} is that superperverse  intersection cohomology modules are independent  of the choice of superperversity below the point at which $\bar p(k)$ becomes $\leq k-2$. In other words, if $\bar p$ and $\bar p'$ are two perversities such that $\bar p(k)$ and $\bar p'(k)$ are both $\geq k-1$ for $k\leq m$ and $\bar p(k)=\bar p'(k)$ for $k>m$, then $I^{\bar p}_{\mf X}H^*(X;\mc G)=I^{\bar p'}_{\mf X}H^*(X;\mc G)$. This is because, as seen in the proof, the condition $\bar p(k)\geq k-1$ for $k\leq m$ implies that $\mc P^*_{m+1}=Ri_*\mc G$, where $i: X-\Sigma\into U_{m+1}$ is the inclusion. We formalize this result as a corollary.

\begin{corollary}\label{C: ultrap}
If $\bar p(k)\geq k-1$ for $k\leq m$, then $\mc P^*_{m+1}=Ri_*\mc G$, where $i: X-\Sigma \into U_{m+1}$ is the inclusion. Thus if $\bar p$ and $\bar p'$ satisfy  $\bar p(k), \bar p'(k)\geq k-1$ for $k\leq m$ and $\bar p(k)=\bar p'(k)$ for $k> m$, then $I^{\bar p}_{\mf X}H^*(X;\mc G)=I^{\bar p'}_{\mf X}H^*(X;\mc G)$. 
\end{corollary}

On the other hand, if $\bar p(k)\geq k-1$ for $k\leq m<n$ but $\bar p(k)\leq k-2$ for $k> m$, then we may obtain intersection cohomology modules that do not come directly either from traditional perversity intersection homology or from ordinary cohomology on complements.  For examples of this phenomenon, see \cite{CS} or \cite{GBF2}.

\begin{remark}
It follows from Corollary \ref{C: ultrap} that we could study all superperversities for intersection cohomology on pseudomanifolds by limiting ourselves to the case $\bar p(1)=0$. However, it may be useful in future notation or when studying more general spaces to retain the more general concept. 
\end{remark}

\section{Where topological invariance breaks down}\label{S: break}

Having seen that superperverse intersection homology is not a topological invariant, we are led to ask where the Goresky-MacPherson  proof of invariance for traditional perversities breaks down. In this section, we will identify the point of difficulty. 

Complete definitions and details of the traditional proof of invariance will be given in the next section. For now we simply note that the idea of the proof is to start with one set of axioms that completely characterize the Deligne sheaf and to progress through a series of equivalent axioms until one reaches a set that does not depend on the choice of stratification. The successive axioms, as treated in Borel \cite{Bo}, are labeled $AX1_{\bar p,\mf X, \mc G}$, $AX1'_{\bar p,\mf X, \mc G}$, $AX2_{\bar p,\mf X, \mc G}$, and $AX2_{\bar p, \mc G}$. Note that each depends on the space $X$, the stratification $\mf X$, the perversity $\bar p$, and the coefficient system $\mc G$, except for the last set of axioms, which does not refer to a specific stratification. 

For a differential graded sheaf $\mc S^*$ on $X$, let $\mc S^*_k=\mc S^*|_{U_k}$. Then $\mc S^*$ satisfies the first set of axioms $AX1_{\bar p,\mf X, \mc G}$ if the following conditions hold:
\begin{description}
\item[$(1a)$] $\mc S^*$ is bounded, $\mc S^i=0$ for $i<0$, and $\mc S^*_1=\mc G$. (Note: since we are really working in the derived category, these conditions can also be stated as  $\mc S^*_1$ is  quasi-isomorphic to $\mc G$ and  $\mc H^i(\mc S^*)=0$ for $i<0$ and $i\gg 0$.)

\item[$(1b)$] For $x\in S_{n-k}$, $k\geq 1$, $H^i(\mc S^*_x)=0$ if $i> \bar p(k)$.

\item[$(1c)$]  The attachment map $\alpha_k:\mc S^*_{k+1}\to Ri_{k*}\mc S^*_k$, $k\geq 1$, is a quasi-isomorphism up to $\bar p(k)$, i.e. it induces isomorphisms on $\mc H^i$ for $i\leq \bar p(k)$. 
\end{description}

We have modified the axioms slightly from their standard form to allow the cases $k=1$ that arise when $X$ has a non-empty pseudoboundary.

It follows as in \cite[Section V.2]{Bo} that any sheaf satisfying these axioms is quasi-isomorphic to the Deligne sheaf with inputs $X$, $\mf X$, $\bar p$, and $\mc G$. This remains true even if $\bar p$ is a superperversity or if $\mf X$ is a stratification with a non-empty pseudoboundary. 
It is also true, even in these more general cases, that the axioms $AX1_{\bar p,\mf X, \mc G}$ are equivalent to the axioms $AX1'_{\bar p,\mf X, \mc G}$ (see Section \ref{S: main}, below). However, for superperversities or stratifications with non-empty pseudoboundaries, these axioms are no longer equivalent to the following axioms $AX2_{\bar p,\mf X, \mc G}$:

\begin{description}
\item[$(2a)$] $\mc S^*$ is bounded, $\mc S^i=0$ for $i<0$, $\mc S^*_1=\mc G$, and $\mc S^*$ is $\mf X$-clc.
\item[$(2b)$] $\dim \{\text{supp} \mc H^j(\mc S^*)\}\leq n-\bar p^{-1}(j)$ for all $j>0$.
\item[$(2c)$] $\dim \{x\in X\mid H^j(f_x^!\mc S^*)\neq 0\}\leq n-\bar q^{-1}(n-j)$ for all $j<n$.
\end{description}

Here $\bar q(k)=k-2-\bar p(k)$ is also a perversity, $\bar p^{-1}(j)=\min\{c\mid \bar p(x)\geq j\}$ if $j\leq \bar p(n)$,  and  $\bar p^{-1}(j)=\infty$  if $j>\bar p(n)$. We similarly define $\bar q^{-1}$.  See Section \ref{S: main} for more details. 

We will demonstrate that, when $\bar p$ is a superperversity, two complexes of sheaves that are not quasi-isomorphic may satisfy $AX2_{\bar p,\mf X, \mc G}$.

Consider the example of the previous section given by the stratified pseudomanifold  $S^n\supset x$, $n\geq 2$. Let $\bar p$ be the superperversity $\bar p(k)=k-1$, and let $\R_{S^n-x}$ be the constant coefficient system on $S^n-x$ with stalk $\R$. The Deligne sheaf for this stratification is $\mc P^*=\tau_{\leq n-1}Ri_*\R_{S^n-x}$, where $i:S^n-x\into S^n$ is the inclusion. As the Deligne sheaf, $\mc P^*$ certainly satisfies $AX1_{\bar p,\mf X, \R_{S^n-x}}$ and hence also $AX1'_{\bar p,\mf X, \R_{S^n-x}}$. 

Let us check that $\mc P^*$ satisfies $AX2_{\bar p,\mf X, \R_{S^n-x}}$.
It is clear that it satisfies $(2a)$. 
Now for $y\neq x \in S^n$,  $\mc H^*(\mc P^*)_y=\R[0]$. At $x$, $\mc H^j(\mc P^*)_x=\lim_{x\in U} \H^j(U-x;\R)\cong H^j(S^{n-1};\R)$, which is $\R$ for $j=0, n-1$, and $0$ otherwise. Thus $\dim \{\text{supp} \mc H^{n-1}(\mc S^*)\}=0$, and $\dim \{\text{supp} \mc H^j(\mc S^*)\}=-\infty$ for $j>0$, $j\neq n-1$. Since $\bar p^{-1}(n-1)=n$, we see that $\mc P^*$ satisfies condition $(2b)$. 

To check $(2c)$, we use the adjunction exact sequence 
\begin{equation}\label{E: les}
\begin{CD}
@>>> & H^j(f_y^!\mc P^*) & @>>>& H^j(\mc P^*_y) & @>\alpha >> & H^j(Ri_{*}i^*\mc P^*)_y &@>>>,
\end{CD}
\end{equation}
where $y$ is any point in $S^n$ and $i: S^n-y\into S^n$. We know from the above computations that for $y\neq x$, $H^*(\mc P^*_y)=\R[0]$, while  $H^j(\mc P^*_x)=\R[0]\oplus \R[-(n-1)]$. Similarly, for $y\neq x$,  $H^j(Ri_{*}i^*\mc P^*)_y= \R[0]\oplus \R[-(n-1)]$. Meanwhile, 
$\alpha: \mc H^j(\mc P^*)_x\cong  \mc  H^j(Ri_{*}i^*\mc P^*)_x$ for $ j\leq \bar p(n)=n-1$ by $(1c)$.  It is also not hard to verify by direct sheaf level computations that the map $\alpha_*: \mc H^0(\mc P^*)_y \to  \mc H^0(Ri_{*}i^*\mc P^*)_y: \R \to \R$ is also an isomorphism  for $y\neq x$. Thus we see that   $H^j(f_y^!\mc P^*)=H^j(f_x^!\mc P^*)=0$ for $j< n$ (though we do have $H^n(f_y^!\mc P^*)=\R$ for $y\neq x$).
So for $j<n$, $\dim \{y\in S^n\mid H^j(f_y^!\mc S^*)\neq 0)=-\infty \leq n-\bar q^{-1}(n-j)$, and $\mc P^*$ satisfies $AX2_{\bar p,\mf X, \R}$

On the other hand, consider the constant sheaf $\R_{S^n}$ on $S^n$ with the same stratification $S^n\supset x$ and superperversity $\bar p$. Since $S^n$ is a manifold, this would be the Deligne sheaf for any traditional perversity. $\R_{S^n}$ also satisfies $AX2_{\bar p,\mf X, \R_{S^n}}$: condition $(2a)$ is clearly satisfied. $H^*(\R_y)=\lim_{y\in U}\H^*(U;\R_U)\cong H^*(\R^n;\R)=\R[0]$ for all $y$ in $S^n$. So for all $j>0$, $\dim \{\text{supp} \mc H^j(\mc S^*)\}=-\infty \leq n-\bar p^{-1}(j)$, and $(2b)$ is satisfied. We also have, by the exact sequence \eqref{E: les} and the computations of the preceding paragraph, that $H^j(f_y^!\R_{S^n}^*)=0$ for $j<n$ and $H^n(f_y^!\R_{S^n}^*)=\R$. So for $j<n$, $\dim \{y\in S^n\mid H^j(f_y^!\R )\neq 0\}=-\infty \leq n-\bar q^{-1}(n-j)$, and $(2c)$ is satisfied. Thus both $\mc  P^*$ and $\R_{S^n}$ satisfy $AX2_{\bar p,\mf X, \R_{S^n-x}}$. However, they are not quasi-isomorphic since $\mc H^{n-1}(\R_{S^n})_x=0$ while $H^{n-1}(\mc P^*)_x=\R$. 

So we must conclude that, for superperversities, $AX2_{\bar p,\mf X, \mc G}$ and $AX1_{\bar p,\mf X, \mc G}$ are not equivalent. In fact, we will see below that $AX1_{\bar p,\mf X, \mc G}$ still implies $AX2_{\bar p,\mf X, \mc G}$, but not conversely. The issue in this example is the following: When $\bar p$ is a traditional perversity, $\bar p(n)\leq n-2$ and so taking $\mc P^*=\tau_{\leq \bar p(n)}Ri_*\R$ in the above example yields a sheaf quasi-isomorphic to $\R$, since the extra cohomology in $\mc H^{n-1}(Ri_*\R)_x$ gets truncated off. However, if $\bar p(n)$ can be  $\geq n-1$, then the $n-1$ dimensional cohomology lives on in $\mc P^*$. As a result, the attachment map $H^{n-1}(\mc P^*_x) \to   H^{n-1}(Ri_{*}i^*\mc P^*)_x$ becomes an isomorphism, instead of the $0$ map that it would be with a traditional perversity, and  
 $H^n(f_x^!\mc P^*)$ becomes $0$, instead of $\R$.  If $\bar p(n)=n-1$, and hence $\bar q(n-1)=-1$, this fact is detected by axiom $(1'c)$, which says that if $x\in S_{n-k}$, $k\geq 1$, then $H^j(f^!_x\mc S^*)=0$ for $j<n-\bar q(k)$. \emph{However}, this subtlety cannot be detected by the condition $(2c)$, since $(2c)$ only looks at $j<n$. This cannot be corrected simply by allowing $j=n$ in $(2c)$ since we have seen that for $y\in S^n-x$, $H^n(f_y^!\mc P^*)=\R$. So we \emph{always} have $\dim \{y\in S^n\mid H^n(f_y^!\mc P^*)\neq 0\}=n$, which  won't be $\leq n-\bar q^{-1}(0)$, since $\bar q^{-1}(0)$ must be positive.  

The study of other examples indicates that the situation described represents the generic difficulty in trying to characterize the Deligne sheaf by axioms of the form of $AX2_{\bar p,\mf X, \mc G}$: when dealing with superperversities, $H^n(f_x^!\mc P^*)$ must somehow be taken into consideration, but $\dim \{y\in S^n\mid H^n(f_y^!\mc P^*)\neq 0\}=n$. The solution, presented in the next section, is to consider only $\dim \{y\in \Sigma \mid H^j(f_y^!\mc S^*)\neq 0\}$. This will rescue the equivalence of $AX2_{\bar p,\mf X, \mc G}$ and $AX1_{\bar p,\mf X, \mc G}$ at the expense of adding an explicit dependence on the choice of $\Sigma$ (though not a dependence on all of $\mf X$).

\begin{remark}
The problem in the case of traditional perversities on stratified pseudomanifolds with nonempty pseudoboundaries is exactly the same. If $X^{n-1}\neq X^{n-2}$, then $\bar p(1)$ comes into play, but if $\bar p(1)=0=1-1$, then we run into the same difficulties as for superperversities in the previous paragraph, those arising from a situation where $\bar p(k)\nleq k-2$. 
\end{remark}

\section{The main theorem}\label{S: main}

We now come to our main theorem.

\begin{theorem}\label{T: main}
Let $X$ be  an $n$-dimensional topological pseudomanifold with (possibly empty)  pseudoboundary. Let $\Sigma$ be the $n-1$  skeleton of some topological stratification of $X$, and let $\mc G$ be a system of local coefficients on $X-\Sigma$. Let $\bar p$ be a traditional perversity or superperversity. Then the Deligne sheaf $\mc P^*\in D^b(X)$ is   independent of the choice of stratification of $X$ subject to $\Sigma$ and hence so are the intersection cohomology modules $I^{\bar p}_{\Sigma}H^*(X;\mc G)$. 
\end{theorem}

In other words, any two topological stratifications  of $X$ such that $X^{n-1}=\Sigma$ yield the same Deligne sheaf $\mc P^*$ up to quasi-isomorphism.  

We shall discuss two proofs.

Our primary exposition basically follows that of Goresky and MacPherson \cite{GM2} for the topological invariance of intersection homology with traditional perversities on topological pseudomanifolds without pseudoboundaries. However, we will follow more closely the treatment of this proof given by Borel in \cite{Bo}, which has the advantage of treating some of the technical issues (particularly sheaf constructibility) slightly more cleanly (see the Remarks in \cite[\S V.4.20]{Bo}). The idea of the proof is to describe the Deligne sheaf by a set of axioms and then progress through several sets of equivalent axioms to one that is no longer dependent upon the specific stratification $\mf X$. As we have seen in the previous sections of this paper, the  equivalence of the usual axioms breaks down when considering superperversities or pseudoboundaries. Our main work then is to modify these axioms in order to reinstate these equivalences. However, we have also seen that it will not be possible to do so in a way that maintains complete stratification independence. Thus we must introduce a dependence on the singular locus $\Sigma$. 

Following this proof, we provide an alternative set of axioms that characterize the Deligne sheaf as a \emph{codimension $\geq c$ intersection cohomology theory} in the sense of Habegger and Saper \cite{HS91}. These axioms have slightly simpler support and cosupport conditions than do the axioms $AX2'_{\bar p,\mf X, \mc G}$ that we present in our first treatment, modified from the Goresky-MacPherson axioms  $AX2_{\bar p,\mf X, \mc G}$. However, the coefficients of the resulting codimension $\geq c$ intersection cohomology theory are more complicated. In fact they are not locally-constant or even clc (though they will be $\mf X$-clc). Consequently, the two axiomatic characterizations we give are of somewhat different character, particularly in their incorporation of the dependence upon the singular locus $\Sigma$. 
 
We begin with the approach following Goresky and MacPherson. We will focus principally on the parts of our proof that diverge from those in \cite{GM2} and \cite{Bo}, though for the purposes of readability and relative completeness, we provide at least an outline of the entire proof. We should also warn the comparing reader that since we are dealing with pseudomanifolds with pseudoboundaries, most of our inductions will start with $k=1$ instead of $k=2$.

N.B. Following the treatment in \cite{Bo}, we suppress much of the derived category notation. However, one should note that an equal sign  between differential graded sheaves, $A^*=B^*$, denotes quasi-isomorphism, i.e. isomorphism in the derived category.

We begin with the first set of axioms $AX1_{\bar p,\mf X, \mc G}$. Here we only require that $\mf X$ be a filtration of $X$ by closed subsets, not necessarily a topological stratification. Let $\mc S^*$ be a differential graded sheaf on $X$, and let $\mc S^*_k=\mc S^*|_{U_k}$, where $U_k=X-X^{n-k}$. Then the axioms $AX1_{\bar p,\mf X, \mc G}$ consist of the following conditions \cite[p. 61, 86]{Bo}:
\begin{description}
\item[$(1a)$] $\mc S^*$ is bounded, $\mc S^i=0$ for $i<0$, and $\mc S^*_1=\mc G$. (Note: since we are really working in the derived category, these conditions can also be stated as  $\mc S^*_1$ is  quasi-isomorphic to $\mc G$ and  $\mc H^i(\mc S^*)=0$ for $i<0$ and $i\gg 0$.)

\item[$(1b)$] For $x\in S_{n-k}$, $k\geq 1$, $H^i(\mc S^*_x)=0$ if $i> \bar p(k)$.

\item[$(1c)$]  The attachment map $\alpha_k:\mc S^*_{k+1}\to Ri_{k*}\mc S^*_k$, $k\geq 1$, is a quasi-isomorphism up to $\bar p(k)$, i.e. it induces isomorphisms on $\mc H^i$ for $i\leq \bar p(k)$. 
\end{description}

The attachment map of condition $(1c)$ is the composition of the natural maps $\mc S_{k+1}^*\to i_{k*}i_k^*\mc S_{k+1}^*\to Ri_{k*}i_k^*\mc S^*_{k+1}=Ri_{k*}\mc S^*_k$. It is automatically a quasi-isomorphism at points $x\in U_k$; the condition implies that for points $x\in U_{k+1}-U_k=S_{n-k}$, then $H^i(\mc S^*_{k+1,x})=\lim_{x\in V}\H^i(V-V\cap S_{n-k};\mc S^*_{k})$ for $i\leq \bar p(k)$, where $V$ is a fundamental set of neighborhoods of $x$ in $U_{k+1}$. 

It is obvious from the definition of the Deligne sheaf $\mc P^*$ that it satisfies this set of axioms. If $\mf X$ is a topological stratification, it is also true that any sheaf $\mc S^*$ satisfying these axioms is quasi-isomorphic to the Deligne sheaf $\mc P^*$. This follows just as in the proof of Theorem V.2.5 of \cite{Bo} via \cite[Lemma V.2.4]{Bo}, which shows that any $\mc S^*$ satisfying these axioms has $\mc S^*_{k+1}=\tau_{\leq \bar p(k)}Ri_{k*}\mc S_k^*$. Furthermore, continuing to assume that $\mf X$ is a topological stratification, it follows from \cite[\S V.3]{Bo} that any such $\mc S^*$ is \emph{$\mf X$-cohomologically constructible} ($\mf X$-cc) and \emph{cohomologically constructible} (cc). In particular, $\mc S^*$ is \emph{$\mf X$-cohomologically locally constant} ($\mf X$-clc).

The next step of the proof is to show the equivalence (under suitable conditions)  of $AX1_{\bar p,\mf X, \mc G}$ with another set of axioms $AX1'_{\bar p,\mf X, \mc G}$. To state these axioms, we need the notion of the dual perversity $\bar q$ such that $\bar p(k)+\bar q(k)=k-2$. Note that so long as $\bar p$ is actually a perversity, i.e. $\bar p(k)\leq \bar p(k+1)\leq \bar p(k)+1$, $\bar q$ will also be a perversity. 
It is  interesting to note that if $\bar p$ were allowed to take jumps of size $>1$ then $\bar q$ might not be monotonic, which would cause difficulties below when we consider $q^{-1}$. For this reason, it is necessary to work with actual perversities and not ``loose perversities'' (see King \cite{Ki}). 

Let $f_x: x\into X$ denote the inclusion of the point $x\in X$, and  let $j_k:S_{n-k}\into U_{k+1}=S_{n-k}\cup U_k$ denote inclusion of the stratum.

$AX1'_{\bar p,\mf X, \mc G}$ reads as follows:
\begin{description}
\item[$(1'a)$] $\mc S^*$ is bounded, $\mc S^i=0$ for $i<0$, $\mc S^*_1=\mc G$, and $\mc S^*$ is $\mf X$-clc. 

\item[$(1'b)$] If $x\in S_{n-k}$, $k\geq 1$, then $H^j(\mc S^*_x)=0$ for $j>\bar p(k)$.

\item[$(1'c)$] If $x\in S_{n-k}$, $k\geq 1$, then $H^j(f^!_x\mc S^*)=0$ for $j<n-\bar q(k)$. 
\end{description}

\begin{remark}
The conditions $k\geq 1$ in the above statements, which would read $k\geq 2$ in the standard case, are left tacit in \cite{Bo} except for the first presentation of $AX1_{\bar p,\mf X, \mc G}$ on page 61. They are certainly assumed, however, if for no other reason than that $\bar p(k)$ is not even defined for $k<2$ in that source. More accurately,  such considerations of $S_n$ are unnecessary anyway, since the first axiom tells us that $\mc S^*|_{U_1}=\mc S^*_1=\mc G$.
We will be careful, however, and include the condition $k\geq 1$ explicitly here, as it will be used below.
\end{remark}

Notice that  $(1'b)$ is the same as $(1b)$, and $(1'a)$ is $(1a)$ together with the condition that $\mc S^*$ be $\mf X$-clc. The following proposition, a slight generalization of  Proposition V.4.3 of \cite{Bo}, continues to hold for superperversities or non-empty pseudoboundaries  (replacing $k\geq 1$ with $k\geq 2$):

\begin{proposition}[Proposition V.4.3 of \cite{Bo}]\label{P: 1 to 1'}
Suppose that each $S_{n-k}$ is a manifold of dimension $n-k$ or is empty, that
$\mc S^*$ is $\mf X$-clc, and that $j_k^!\mc S^*$ is cohomologically locally constant (clc) for $1\leq k\leq n$. Then $\mc S^*$ satisfies $AX1_{\bar p,\mf X, \mc G}$ if and only if it satisfies $AX1'_{\bar p,\mf X, \mc G}$.
\end{proposition}

The proof of this proposition given in \cite{Bo} continues to hold despite our generalizations. The key point is the equivalence of $(1c)$ and $(1'c)$ in the presence of $(1a)$, $(1b)=(1'b)$, and the additional hypotheses of the proposition. This equivalence is demonstrated by showing that $(1c)$ and $(1'c)$ are each equivalence to $(1''c)$:

\begin{description}
\item[$(1''c)$] If $x\in S_{n-k}$, then $H^j(j^!_k\mc S^*)=0$ for $j\leq \bar p(k)+1$. 
\end{description}

The equivalence of $(1c)$ and $(1''c)$, given $(1b)$, follows from the exact sequence
\begin{equation*}
\begin{CD}
@>>> & \mc H^j(j_k^!\mc S^*)_x & @>>>&\mc  H^j(\mc S^*)_x & @>\alpha^j_k >> & \mc H^j(Ri_{k*}\mc S^*_k)_x &@>>>
\end{CD}
\end{equation*}
for $x\in S_{n-k}$.
The equivalence of $(1'c)$ and $(1''c)$ uses the equality $H^j(f_x^!\mc S^*)=H^j(\ell_x^!(j_k^!\mc S^*))=H^{j-n+k}((j_k^!\mc S^*)_x)$, where $\ell_x: x\into S_{n-k}$ is the inclusion and the second equality is due to \cite[Proposition V.3.7.b]{Bo}, owing to  $j_k^!\mc S^*$ being clc. 

It follows from Proposition \ref{P: 1 to 1'} that if $X$ (with filtration $\mf X$) is an $n$-dimensional stratified topological pseudomanifold with pseudoboundary, then $AX1_{\bar p,\mf X, \mc G}$ and $AX1'_{\bar p,\mf X, \mc G}$ are equivalent: If $\mf X$ is a topological stratification then each $S_{n-k}$ is a manifold of dimension $n-k$ (or empty). If $\mc S^*$ satisfies $AX1_{\bar p,\mf X, \mc G}$, then it is $\mf X$-clc since it is quasi-isomorphic to the Deligne sheaf, while being $\mf X$-clc is an explicit condition of $(1'a)$. Then by \cite[Proposition V.3.10.d]{Bo}, if $\mc S^*$ is $\mf X$-clc then $j^!_k\mc S^*$ is clc.  Proposition V.3.10.d of \cite{Bo} is stated for stratified pseudomanifolds, but it holds equally well for pseudomanifolds with pseudoboundary. 

Thus if $\mf X$ is a topological stratification, any sheaf $\mc S^*$ on $X$ satisfying $AX1'_{\bar p,\mf X, \mc G}$ is quasi-isomorphic to the Deligne sheaf.

\medskip

The next stage of the program is the one at which we must begin to make some changes, as we have already seen  that, for superperversities,  
$AX1_{\bar p,\mf X, \mc G}$ is not equivalent to $AX2_{\bar p,\mf X, \mc G}$. To state these next axioms, we need the concept of $\bar p^{-1}$. To simplify somewhat our working with these inverses, let us extend the domain of the perversity $\bar p$ to all of $\Z$. This extension is simply for notational convenience so that we do not need to be as careful with the input numbers to our perversities. Suppose that $\bar p$ is already defined on $\Z^{\geq 1}$. Then for $k<1$, we let $\bar p(k)=\bar p(1)+k-1$, and for $k>n$, let $\bar p(k)=\bar p(n)$. If we are in the traditional case where $X^{n-1}=X^{n-2}$ and our given $\bar p$ is only defined for $k\geq 2$, then take $\bar p(1)=\bar p(2)-1$ and then define $\bar p$ on $\Z$ as in the previous case.  $\bar q$ continues to be defined by $\bar p(k)+\bar q(k)=k-2$ and so its domain is also extended to all of $\Z$. 
 
Now we can define $\bar p^{-1}$ on all of $\Z$ (cf. \cite[p. 88]{Bo}).
If $j\leq \bar p(n)$, let $\bar p^{-1}(j)=\min\{c\mid \bar p(x)\geq j\}$. If $j>\bar p(n)$, let $\bar p^{-1}(j)=\infty$. $\bar q^{-1}$ is defined similarly. Then for all $k\in \Z$, 
\begin{equation}
\bar p(k)\geq j \Leftrightarrow k\geq \bar p^{-1}(j).
\end{equation}
The following two useful formulas follow immediately:
\begin{equation}\label{E: p}
j \leq \bar p(k) \Leftrightarrow n-k\leq  n- \bar p^{-1}(j)
\end{equation}
\begin{equation}\label{E: q}
j \geq n-\bar q(k)  \Leftrightarrow n-k\leq  n- \bar q^{-1}(n-j).
\end{equation}
Here we take $\bar q^{-1}(j)=-\infty$ if $j\leq \bar q(k)$ for all $k$ (or, equivalently, if $j\leq \bar q(1)$).

The standard set of axioms  $AX2_{\bar p,\mf X, \mc G}$ read as follows:
\begin{description}
\item[$(2a)$] $\mc S^*$ is bounded, $\mc S^i=0$ for $i<0$, $\mc S^*_1=\mc G$, and $\mc S^*$ is $\mf X$-clc.
\item[$(2b)$] $\dim \{\text{supp} \mc H^j(\mc S^*)\}\leq n-\bar p^{-1}(j)$ for all $j>0$.
\item[$(2c)$] $\dim \{x\in X\mid H^j(f_x^!\mc S^*)\neq 0\}\leq n-\bar q^{-1}(n-j)$ for all $j<n$.
\end{description}

For the intersection cohomology with traditional perversities of  topological pseudomanifolds without pseudoboundaries, these axioms are equivalent to  $AX1'_{\bar p,\mf X, \mc G}$ and hence also to $AX1_{\bar p,\mf X, \mc G}$.
As we have seen in Section \ref{S: break}, this is not the  case for  superperversities or pseudomanifolds with pseudoboundary. The difficulty is the restriction $j<n$ in condition (2c). In the standard proof of equivalence (see \cite[p. 89]{Bo} or our modification below), it is shown  that $(2c)$ and $(1'c)$ are equivalent under the appropriate hypotheses. The implication $(1'c)\Rightarrow (2c)$ continues to hold in the more general case, but the reverse implication becomes insufficient since it is no longer enough to consider what happens only when $j<n$. 
For a \emph{traditional} perversity $\bar p$ on a topological pseudomanifold without pseudoboundary, $S_{n-k}=\emptyset$ for $k\notin [2,n]$, and, for $k$ in this range, $\bar p(k)$ and $\bar q(k)$ are both $\geq 0$. Thus the condition $j<n-\bar q(k)$ in $(1'c)$ implies that $j<n$, allowing the restriction to this range in $(2c)$. Conversely, this range of parameters in $(2c)$ covers all possibilities in $(1'c)$.  However, if $\bar p$ is a superperversity or $X$ has a non-empty pseudoboundary, we might have $\bar q(k)<0$ for some relevant $k\geq 1$, in which case $n-\bar q(k)\geq n$. Thus in order to obtain condition $(1'c)$ from $(2c)$, it is necessary to modify $(2c)$ to take into account the cases $j\geq n$.

Unfortunately, simply removing the condition $j<n$ in $(2c)$ is insufficient. For $j=n$, $\dim \{x\in X\mid H^n(f_x^!\mc S^*)\neq 0\}=n$, since 
if $G$ is the stalk of $\mc G$, $H^n(f_x^!\mc S^*)=G$ for $x$ in the $n$-manifold $U_1$.  If we let $i:X-x\into X$ be the inclusion, this follows from the long exact sequence 
\begin{equation*}
\begin{CD}
@>>> & H^j(f_x^!\mc S^*) & @>>>& H^j(\mc S^*_x) & @>\alpha >> & H^j((Ri_{*}i^*\mc S^*)_x)&@>>>
\end{CD}
\end{equation*}
since we know that  $H^j(\mc S^*_x)=0$ for $j>0$, $x\in U_1$, and  $H^j((Ri_{k*}\mc S^*_k)_x)=\lim_{x\in U}\H^j(U-x; \mc G)=H^j(\mc S^{n-1}; G)$. Thus   $H^n(f_x^!\mc S^*)=G$ if $n\geq 2$. 

So we cannot hope to satisfy any analogue of $(2c)$ that  would include a conditions  $\dim \{x\in X\mid H^n(f_x^!\mc S^*)\neq 0\}\leq n-\bar q^{-1}(0)$, at least not if  $\bar q^{-1}(0)>0$, which will always hold if $\bar p$ is a superperversity. Thus we modify $(2c)$ in such a way to exclude this problematic dense set $U_1$ by only looking in $X^{n-1}$. The debt we pay is that including  $X^{n-1}$ in this manner, we will not be able to remove it from other axioms later. 

We replace  $AX2_{\bar p,\mf X, \mc G}$ with the following modified axioms $AX2'_{\bar p,\mf X, \mc G}$:
\begin{description}
\item[$(2'a)$] $\mc S^*$ is bounded, $\mc S^i=0$ for $i<0$, $\mc S^*_1=\mc G$, and $\mc S^*$ is $\mf X$-clc.
\item[$(2'b)$] $\dim \{\text{supp} \mc H^j(\mc S^*)\}\leq n-\bar p^{-1}(j)$ for all $j>0$.
\item[$(2'c)$] for all $j\neq n$, $\dim \{x\in X\mid H^j(f_x^!\mc S^*)\neq 0\}\leq n-\bar q^{-1}(n-j)$ ; $\dim \{x\in X^{n-1}\mid H^n(f_x^!\mc S^*)\neq 0)\leq n-\bar q^{-1}(0)$.
\end{description}

We can now prove the following:

\begin{proposition}\label{P: AX1 to Ax2'}
Suppose that $\mf X$ is such that each $S_{n-k}$ is a manifold of dimension $n-k$ or empty and that $\mc S^*$ is a differential graded sheaf such that $j^!_k\mc S^*$ is clc for each $k\geq 1$. Then $\mc S^*$ satisfies $AX2'_{\bar p,\mf X, \mc G}$ if and only if it satisfies $AX1'_{\bar p,\mf X, \mc G}$. In particular, if $\mf X$ is a topological stratification, then $AX2'_{\bar p,\mf X, \mc G}$ and  $AX1'_{\bar p,\mf X, \mc G}$ are equivalent and characterize $\mc S^*$ uniquely up to quasi-isomorphism as the Deligne sheaf. 
\end{proposition}

As this proposition contains the heart of our modification to the standard treatment, we provide a full proof. However, the proof is almost identical to that of \cite[Proposition V.4.9]{Bo}. Our main alterations are the addition of the cases $j=n$ in $(2'c)\Rightarrow (1'c)$ and, of course, the addition of the cases $k=1$. 

\begin{proof}
Since $(1'a)=(2'a)$, we need only show that $(1'b)\Leftrightarrow (2'b)$ and $(1'c)\Leftrightarrow(2'c)$, given   the other hypotheses . 

\paragraph{$(1'b)\Rightarrow (2'b)$:} For $x\in S_n=U_1$, $\mc S^*_1=\mc G$ by $(1'a)$, so $H^*(\mc S^*_x)=G[0]$. If $x\in S_{n-k}\neq \emptyset$, $k\geq 1$, and $H^j(\mc S_x^*)\neq 0$, then $j\leq \bar p(k)$ by $(1'b)$. But if $j\leq \bar p(k)$ then $n-k\leq n-\bar p^{-1}(j)$ by equation \eqref{E: p}. Since $S_{n-k}\neq \emptyset$, $\dim S_{n-k}=n-k$ by hypothesis, and So $\dim (\text{supp }\mc H^j(\mc S^*))\leq n-\bar p^{-1}(j)$ for $j>0$. 

 \paragraph{$(2'b)\Rightarrow (1'b)$:} Note that $\bar p(k)\geq 0$, so to see that $H^j(\mc S^*_x)=0$ for $j>\bar p(k)$, it suffices to consider $j>0$. If $x\in S_{n-k}$, $k\geq 1$,  and $H^j(\mc S^*_x)\neq 0$, $j>0$, then, since $\mc S^*$ is $\mf X$-clc, $H^j(\mc S^*_y)\neq 0$ for all $y$ in the same $n-k$ dimensional connected component of $S_{n-k}$ as $x$. Thus by $(2'b)$ we must have $n-k\leq n-\bar p^{-1}(j)$, which implies by equation \eqref{E: p} that $j\leq \bar p(k)$. So if $j>\bar p(k)\geq 0$, $H^j(\mc S^*_x)=0$.

\paragraph{$(1'c)\Rightarrow (2'c)$:} For $x\in S_n=U_1$, $\mc S^*_1=\mc G$ by $(1'a)$, so $H^*(f_x^!\mc S^*)=G[-n]$. If $x\in S_{n-k}\neq \emptyset$, $k\geq 1$, and $H^j(f^!_x\mc S^*)\neq 0$, then $j\geq n-\bar q(k)$ by $(1'c)$. But if $j\geq n-\bar q(k)$, then $n-k\leq n-\bar q^{-1}(n-j)$ by equation \eqref{E: q}. Since $S_{n-k}\neq \emptyset$,  $\dim S_{n-k}=n-k$ by hypothesis, and  so for $j\neq n$, $\dim\{x\in X\mid  H^j(f^!_x \mc S^*)\neq 0\} \leq n-\bar q^{-1}(n-j)$. If $j=n$, $\dim H^n(f_x^!\mc S^*)=n$, but the above argument still implies $\dim\{x\in X^{n-1}\mid  H^j(f^!_x \mc S^*)\neq 0\} \leq n-\bar q^{-1}(0)$ since $S_{n-k}\subset X^{n-1}$ for $k\geq 1$.   

\paragraph{$(2'c)\Rightarrow (1'c)$:} Suppose $x\in S_{n-k}$, $k\geq 1$, and that $\ell_x:x\into S_{n-k}$ is the inclusion. Then $f_x^!=\ell_x^!j_k^!$. Since $j_k^!\mc S^*$ is clc by hypothesis, we have $H^j((j_k^!\mc S^*)_x)= H^j((j_k^!\mc S^*)_y)$ for all $j$ and for all $y$ in the same connected component of $S_{n-k}$ as $x$. Furthermore, by \cite[Proposition V.3.7.b]{Bo}, since $S_{n-k}$ is a manifold of dimension $n-k$, $H^{j-n+k}((j_k^!\mc S^*)_x)= H^j(\ell_x^!j_k^!\mc S^*)= H^j(f_x^!\mc S^*)$. Thus also $H^j(f_x^!\mc S^*)= H^j(f_y^!\mc S^*)$ for all $y$ in the same connected component of $S_{n-k}$ as $x$. In particular, if $H^j(f_x^!\mc S^*)\neq 0$ and $j\neq n$, then by $(2'c)$, $n-k\leq n-\bar q^{-1}(n-j)$, which implies by equation \eqref{E: q} that  $j\geq n-\bar q(k)$. So for $j\neq n$ and $k\geq 1$, if $x\in S_{n-k}$, $H^j(f_x^!\mc S^*)=0$ for $j<n-\bar q(k)$. 

It remains to consider $j= n$. We already know that $\dim\{x\in X\mid  H^n(f^!_x \mc S^*)\neq 0\}=n$ since $\mc S^*_1=\mc G$ by $(2'a)$. But $(1'c)$ only concerns the case  $k\geq 1$, i.e. the strata $S_{n-k}$ that are contained in $X^{n-1}$. In $X^{n-1}$, the second clause of $(2'c)$ says that  $\dim \{x\in X^{n-1}\mid H^n(f_x^!\mc S^*)\neq 0\}\leq n-\bar q^{-1}(0)$. So by the same argument as in the previous paragraph, if $x\in S_{n-k}$ and $H^n(f_x^!\mc S^*)\neq 0$, then $n-k\leq n-\bar q^{-1}(0)$and $n\geq n-\bar q(k)$. So if $x\in S_{n-k}$, $k\geq 1$, then $H^n((f_x^!\mc S^*)=0$ for $n<n-\bar q(k)$ (i.e. when $\bar q(k)<0$). 

Thus $(1'c)$ holds for all $j$. 

\medskip

Finally, if $\mf X$ is a topological stratification, then each $S_{n-k}$ is indeed a manifold of dimension $n-k$ or empty and, since $\mc S^*$ must be $\mf X$-clc by either set of axioms, each $j_k^!\mc S^*$ is clc for all $k$ by \cite[Proposition V.3.10.d]{Bo}. Thus $AX2'_{\bar p,\mf X, \mc G}\Leftrightarrow AX1'_{\bar p,\mf X, \mc G}\Leftrightarrow AX1_{\bar p,\mf X, \mc G}$, which characterizes $\mc S^*$ as the Deligne sheaf.
\end{proof}

The next step in  Borel's treatment of topological invariance in \cite{Bo} is to consider possible coefficient systems $\mc G$ and stratifications of $X$ that are adapted to $\mc G$. Since such adaptation issues concern only the dense top stratum of a stratified pseudomanifold, which in our case will be determined \emph{a priori} by our choice of $X^{n-1}$, we can avoid this discussion. Once we have chosen a top stratum $X^{n-1}=\Sigma$, we are forced to work with coefficient systems $\mc G$ which are defined on (or can be extended uniquely to) $X-\Sigma$. 

This brings us to the axioms $AX2''_{\bar p,\Sigma, \mc G}$. We continue to assume that $X$ is  an $n$-dimensional topological pseudomanifold with (possibly empty)  pseudoboundary. Let $\Sigma$ be the $n-1$  skeleton of some topological stratification of $X$, and let $\mc G$ be a system of local coefficients on $X-\Sigma$. We continue to let  $\bar p$ be a traditional or super-perversity. Then we define the axioms  $AX2''_{\bar p,\Sigma, \mc G}$ on a differential graded sheaf $\mc S^*$ on $X$ as follows:
\begin{description}
\item[$(2''a)$] $\mc S^*$ is bounded, $\mc S^i=0$ for $i<0$, $\mc S^*|_{X-\Sigma}=\mc G$, and $\mc S^*$ is $\mf X$-clc for some topological stratification of $X$ subject to $\Sigma$.

\item[$(2''b)$] $\dim \{\text{supp} \mc H^j(\mc S^*)\}\leq n-\bar p^{-1}(j)$ for all $j>0$.

\item[$(2''c)$] for all $j\neq n$, $\dim \{x\in X\mid H^j(f_x^!\mc S^*)\neq 0\}\leq n-\bar q^{-1}(n-j)$ ;  $\dim \{x\in \Sigma\mid H^n(f_x^!\mc S^*)\neq 0\}\leq n-\bar q^{-1}(0)$.
\end{description}

This is our analogue to the axioms $AX2_{\bar p, \mc G}$ in the standard treatments (see \cite{GM2}, \cite[p. 90-91]{Bo}). Notice that the axioms $AX2''_{\bar p,\Sigma, \mc G}$ do not depend on any particular $\mf X$, but they do depend on $\Sigma$. The latter two axioms are essentially the same as those in $AX2'_{\bar p,\mf X, \mc G}$.

The following is our analogue of \cite[Theorem 4.15]{Bo}:

\begin{theorem}\label{T: unterstrat}
Let $X$, $\Sigma$, $\bar p$, and $\mc G$ be as above. Then there exists a differential graded sheaf $\td{\mc P}^*$ satisfying $AX2''_{\bar p,\Sigma, \mc G}$ and $AX2'_{\bar p,\mf X, \mc G}$ for every topological stratification $\mf X$ subject to $\Sigma$. 
\end{theorem}
\begin{proof}
Once again, the proof consists mostly of minor modifications to that presented for \cite[Theorem 4.15]{Bo}. We provide an outline of the entire proof, noting our deviations in somewhat more detail.

$\td{\mc P}^*$ will be the Deligne sheaf associated to $\mc G$, $\bar p$,  and a filtration $\td{\mf X}$ $$X=\td X^n\supset \td X^{n-1}=\Sigma \supseteq \cdots\supseteq \td X^0\supset \td X^{-1}=\emptyset$$
that satisfies the properties $(I_k)$ and $(II_k)$ listed below. As for our previous notation, let  $\td U_k=X-\td X^{n-k}$, $\td S_{n-k}=\td X^{n-k}-\td X^{n-k-1}$, $\td i_k:\td U_k\into \td U_{k+1}$, $\td j_k:\td S_{n-k}\into \td U_{k+1}$, $\td{\mc P}^*_1=\mc G$ on $\td U_1$, $\td{\mc P}^*_{k+1}=\tau_{\leq \bar p(k)}R\td i_{k*}\td{\mc P}_k^*$ for $k\geq 1$, and $\td{\mc P}^*=\td{\mc P}^*_{n+1}$. Also note that $\td X^{n-1}=\Sigma$. Then for each $k$, $0\leq k\leq n$, we want
\begin{description}

\item[$(I_k)$] \begin{enumerate}

\item $\td S_{n-k}$ is a manifold of dimension $n-k$ or is empty.

\item $\td j_k^*\td{\mc P}^*_{k+1}$ is clc. 

\item $\td j_k^!\td{\mc P}^*_{k+1}$ is clc.

\end{enumerate}

\item [$(II_k)$] For every topological stratification $\mf X$ of $X$ subject to $\Sigma$, $\td S_{n-k}$ is a union of connected components of strata of $\mf X$ and $U_{k+1}\subset \td U_{k+1}$.  
\end{description}

If such an $\td{\mf X}$ exists, then the proof of the theorem concludes as follows: $\td{\mc P}^*$ satisfies $AX1_{\bar p,\td{\mf X}, \mc G}$ by construction. So by the above equivalences of axioms  and by the condition $(I_k)$ for all $k$, $1\leq k\leq n$, $\td{\mc P}^*$ also satisfies $AX2'_{\bar p,\td{\mf X}, \mc G}$. Suppose now that $\mf X$ is a topological stratification of $X$ subject to $\Sigma$. Such an $\mf X$ exists since the existence of $\Sigma$ is predicated upon it. Note that $\Sigma=\td X^{n-1}=X^{n-1}$. By conditions $(I)$ and $(II)$, it follows that $\td{\mc P}^*$ is also $\mf X$-clc. Therefore $\td{\mc P}^*$ also satisfies $AX2'_{\bar p,\mf X, \mc G}$. Since $\mf X$ was an arbitrary topological stratification subject to $\Sigma$, it follows that $\td{\mc P}^*$ satisfies $AX2''_{\bar p,\Sigma, \mc G}$ and $AX2'_{\bar p,\mf X, \mc G}$ for all $\mf X$ subject to $\Sigma$. 

The remainder of the proof of the theorem concerns the construction of an $\td{\mf X}$ satisfying properties $(I)$ and $(II)$. This is done precisely as in \cite[pp. 92-93]{Bo}, given that $\td X^{n-1}=\Sigma$. We summarize: 

Since we need only consider topological stratifications $\mf X$ subject to $\Sigma$, $\td U_1=X-\Sigma$ will always be the union of connected components of the  $n$-dimensional strata of $X-X^{n-1}$ for any such $\mf X$. In fact, for any $\mf X$ subject to $\Sigma,$ $ U_1=\td U_1$ and $S_n=\td S_n$. It is clear then that $(I_0)$ and $(II_0)$ hold.

Assume now by induction that $\td U_i$ has been defined for $0<i\leq k$, $k>0$, that $\td X^{n-i}=X-\td U_i$ for $0<i\leq k$, that $\td S_{n-k}=\td X^{n-i}-\td X^{n-i-1}$ for $0\leq i<k$,  and that $(I_i)$ and $(II_i)$ are satisfied for $0\leq i<k$. If we now choose $\td S_{n-k}$, then  we can let  $\td U_{k+1}=\td U_k \cup \td S_{n-k}$. So, we must define $\td S_{n-k}$ and show that $(I_k)$ and $(II_k)$ hold. 

Since we have $\td U_i$ defined  for $i\leq k$, $\td{\mc P}^*_k$ is defined. If $\bar i_k:\td U_k\to X$ and $\bar j_k:\td X^{n-k}\to X$ are the inclusions, define $\bar{\mc P}^*_{k+1}=\tau_{\leq \bar p(k)}R\bar i_{k*}\td{\mc P}^*_k$. Let $\td S'_{n-k}$ be the largest submanifold of $\td X^{n-k}$ of dimension $n-k$, let $\td S_{n-k}''$ be the largest open subset of $\td X^{n-k}$ over which $\bar j_k^*\bar{\mc P}^*_{k+1}$ is clc, and let $\td S_{n-k}'''$ be the largest open subset of $\td X^{n-k}$ over which $\bar j_k^!\bar{\mc P}^*_{k+1}$ is clc. Then we take $\td S_{n-k}=\td S_{n-k}'\cap \td S_{n-k}'' \cap \td S_{n-k}'''$, $\td U_{k+1}=\td U_k\cup \td S_{n-k}$, and $\td{\mc P}^*_{k+1}=\tau_{\leq \bar p(k)}R \td i_{k*}\td{\mc P}^*_k=\bar{\mc P}^*_{k+1}|_{\td U_{k+1}}$. It is clear from these choices and induction that $\td U_{k+1}$ is open and that $(I_k)$ holds at this stage. 

The condition $(II_k)$ follows exactly as in the proof of \cite[Lemma V.4.16]{Bo}, which states that if $Y$ is a connected component of a stratum of some topological stratification $\mf X$ then the intersection of $Y$ with each of $\td S_{n-k}'$, $\td S_{n-k}''$, and $ \td S_{n-k}'''$ is either empty or all of $Y$ and that if $Y\not\subset \td U_k$ and $\text{codim}_X Y=k$ then $Y\subset \td S_{n-k}'\cap \td S_{n-k}'' \cap \td S_{n-k}'''$. This lemma is proven by working locally in distinguished neighborhoods to show that if $y\in Y$, then $y$ has a neighborhood in $Y$ that is either contained completely in $\td S_{n-k}'$ or is disjoint from it and similarly for $\td S_{n-k}''$ and $ \td S_{n-k}'''$. See \cite[pp. 93-4]{Bo}. 
\end{proof}

The rest of the proof of Theorem \ref{T: main} follows immediately: By Theorem \ref{T: unterstrat}, given  $X$, $\Sigma$, $\bar p$, and $\mc G$, there exists a differential graded sheaf $\td P^*$ satisfying $AX2''_{\bar p,\Sigma, \mc G}$ and $AX2'_{\bar p,\mf X, \mc G}$ for every topological stratification $\mf X$ subject to $\Sigma$. But we have seen that if a differential graded sheaf satisfies $AX2'_{\bar p,\mf X, \mc G}$, then it satisfies $AX1_{\bar p,\mf X, \mc G}$, which implies that it quasi-isomorphic to the Deligne sheaf associated to the input data. In other words, $\td{\mc P}^*$ is quasi-isomorphic to all of the Deligne sheaves $\mc P^*$ over all possible $\mf X$ subject to $\Sigma$ ($\bar p$ and $\mc G$ also being fixed). Thus these sheaf complexes are all quasi-isomorphic to each other. Since intersection cohomology of a stratified topological pseudomanifold is the hypercohomology of the Deligne sheaf, it follows that intersection cohomology is independent of choice of stratification $\mf X$ of $X$ subject to $\Sigma$, and $I^{\bar p}_{\Sigma}H^*(X;\mc G)$ is well-defined. \qedsymbol

\paragraph{An alternative characterization.}\label{S: codimc} 

Here we present an alternative set of axioms. These allow us to provide a second conclusion to  the proof of Theorem \ref{T: main} by invoking the codimension $\geq c$ intersection cohomology theory of Habegger and Saper \cite{HS91}. We will call these axioms $AX3_{\bar p,\mf X, \mc G}$, though they should not be confused with the axioms $AX3$ of Goresky and MacPherson \cite{GM2}, which characterize traditional perversity intersection homology in terms of its duality properties. 

Let $\bar p$ be a fixed perversity. Let $c_{\bar p}$, or simply $c$, denote $\bar q^{-1}(0)$, which is $\min \{k\in \Z\mid \bar p(k)\leq k-2\}$ or $\infty$ if $\bar p(k)>k-2$ for all $k$. Then we let $AX3_{\bar p,\mf X, \mc G}$ be  the following set of axioms:

\begin{description}
\item[$(3a)$] $\mc S^*$ is bounded, $\mc S^i=0$ for $i<0$, $\mc S^*$ is $\mf X$-clc, and $\mc S^*_c=(Ri_*\mc G)|_{U_c}$, where $i:X-\Sigma\into X$ is the inclusion.
\item[$(3b)$] $\dim \{\text{supp} \mc H^j(\mc S^*)\}\leq n-\bar p^{-1}(j)$ for all $j>c-2$.
\item[$(3c)$] for all $j<n$, $\dim \{x\in X\mid H^j(f_x^!\mc S^*)\neq 0\}\leq n-\bar q^{-1}(n-j)$.
\end{description}

In the language of \cite{HS91}, these conditions say that $\mc S^*$ is a codimension $\geq c$ intersection cohomology theory with coefficients in $Ri_*\mc G$. We will say more about this below; first we show that on a stratified topological pseudomanifold with (possibly empty) pseudoboundary, these axioms are equivalent to those already studied.

\begin{proposition}\label{P: AX2' to AX3}\label{P: 2' to 3}
If $\mf X$ is a topological stratification of $X$, then $AX3_{\bar p,\mf X, \mc G}$ and  $AX2'_{\bar p,\mf X, \mc G}$ are equivalent. Hence  $AX3_{\bar p,\mf X, \mc G}$ characterizes $\mc S^*$ uniquely up to quasi-isomorphism as the Deligne sheaf. 
\end{proposition}

\begin{proof}
We first show that $AX2'_{\bar p,\mf X, \mc G}$ implies $AX3_{\bar p,\mf X, \mc G}$. 

Since $c_{\bar p}$ must be $\geq 2$ for any traditional perversity or superperversity, $(2'b)$ and $(2'c)$ immediately imply $(3b)$ and $(3c)$. Also, $(2'a)$ implies all of $(3a)$ except for the statement $\mc S^*_c=Ri_*\mc G|_{U_c}$. However, if $\mc S^*$ satisfies $AX2'_{\bar p,\mf X, \mc G}$ for a topological stratification $\mf X$, we already know that $\mc S^*$ is quasi-isomorphic to the Deligne sheaf. Thus $\mc S^*_c=Ri_*\mc G|_{U_c}$, from the proof of Proposition \ref{P: ultrap}. Note that the value $m$ in the statement of Proposition \ref{P: ultrap} must be less than $c$.

Next we show the reverse implication. It is clear that $(3a)$ implies $(2'a)$, that $(3b)$ implies $(2'b)$ for $j>c-2$, and that $(3c)$ implies $(2'c)$ for $j<n$. 

Next we show that $AX3_{\bar p,\mf X, \mc G}$ implies $(2'b)$ for $j\leq c-2$. So suppose that $j\leq c-2$. Then  $\bar p(j+1)\geq j$, so $j+1 \geq \bar p^{-1}(j)$ and $n-j-1\leq n-\bar p^{-1}(j)$. Thus it suffices to show that for $j\leq c-2$, $\dim \{\text{supp} \mc H^j(\mc S^*)\}\leq n-j-1$. We first observe that since $j\leq c-2$, $n-c\leq n-j-2$. So $\mc H^j(\mc S^*)$ can take any values on $X^{n-c}$ without $\dim \{\text{supp} \mc H^j(\mc S^*)\}\leq n-j-1$ being violated, and it suffices to show that the dimension of the intersection of the support of $\mc H^j(\mc S^*)$ with $U_c$ is  $\leq n-\bar p^{-1}(j)$.
For this, condition $(3a)$ tells us that $\mc S^*|_{U_c}=Ri_*\mc G|_{U_c}$. So if we restrict attention completely to the pseudomanifold $U_c$ stratified by the restriction of $\mf X$, then $\bar p(k)\geq k-1$ for all strata $S_{n-k}$ in $U_c$, and, by the proof of Proposition \ref{P: ultrap}, $Ri_*\mc G$ is the Deligne sheaf on $U_c$ (where we restrict $i$ to $i: U_c-\Sigma\into U_c$). Thus, in particular, $AX2'_{\bar p,\mf X|_{U_c}, \mc G}$ holds on $U_c$, which implies that $\dim \{\text{supp} \mc H^j(\mc S^*)\}\cap U_c \leq n-\bar p^{-1}(j)$, as desired.

Finally, we must show that $\dim \{x\in X\mid H^j(f_x^!\mc S^*)\neq 0\}\leq n-\bar q^{-1}(n-j)$ for $j>n$ and  $\dim \{x\in X^{n-1}\mid H^n(f_x^!\mc S^*)\neq 0\}\leq n-\bar q^{-1}(0)$. The first part holds immediately since $\bar q^{-1}(k)=-\infty$ when $k<0$. For the latter statement, recall that $\bar q^{-1}(0)=c$, so $n-\bar q^{-1}(0)=n-c$. It follows as for the last condition that it suffices to show that $\dim \{x\in X^{n-1}\cap U_c\mid H^n(f_x^!\mc S^*)\neq 0\}\leq n-c$. But again this holds since $\mc S^*|_{U_c}=(Ri_*\mc G)_{U_c}$, which is the Deligne sheaf on $U_c$. So 
$AX2'_{\bar p,\mf X|_{U_c}, \mc G}$ holds, which implies by $(2'c)$ that $\dim \{x\in X^{n-1}\cap U_c\mid H^n(f_x^!\mc S^*)\neq 0\}\leq n-c$.

The final claim, that   $AX3_{\bar p,\mf X, \mc G}$ characterizes $\mc S^*$ uniquely up to quasi-isomorphism as the Deligne sheaf, is now a consequence of Proposition \ref{P: AX1 to Ax2'}.
\end{proof}

Since Proposition \ref{P: 2' to 3} implies  that the Deligne sheaf satisfies $AX3_{\bar p,\mf X, \mc G}$, it follows that the Deligne sheaf is a codimension $\geq c$ intersection cohomology theory  with coefficients in $Ri_*\mc G$, where $i:X-\Sigma \into X$.  We will not give the most general definition of  a codimension $\geq c_{\bar p}$ intersection cohomology theory (c-ICT, for short); we refer the reader to \cite{HS91} for complete details. However, we do observe that condition $(3b)$ and $(3c)$ constitute the theory denoted $A_{\bar p}$ in \cite{HS91}, and it is shown there that $A_{\bar p}$ satisfies the more general axioms of  a $c$-ICT. 
By Definition  6.1 of \cite{HS91}, the $A_{\bar p}$ intersection cohomology sheaf with coefficients in $\mc E$ is the unique (up to canonical isomorphism in $D^b(R_X)$) $\mf X$-cc extension of $\mc E|_{U_c}$  to all of $X$ that satisfies $A_{\bar p}$ (the uniqueness follows from the fact that $A_{\bar p}$ is a c-ICT and by the definition of a c-ICT; see \cite[pp. 255-258]{HS91}). Since the Deligne sheaf satisfies $A_{\bar p}$, it must be the unique $\mf X$-cc extension of $\mc P^*|_{U_c}$ that does so. Therefore, since $\mc P^*|_{U_c}=(Ri_*\mc G)|_{U_c}$, the coefficients of $\mc P^*$ are $Ri_*\mc G$ (or any other $\mf X$-cc extension of $(Ri_*\mc G)_{U_c}$ to $X$).  Note that this coefficients system depends on the choice of $\Sigma$, but it does not depend on further refinements of the stratification. 

The conclusion of the proof of our Theorem \ref{T: main} is then also a consequence of  Theorem 6.2 of \cite{HS91}, which states that  the $A_{\bar p}$ intersection cohomology sheaf with coefficients in $\mc E$ is independent of stratifications of $X$ adapted to $\mc E$. The requirement of a stratification $\mf X$ being adapted to $\mc E$ means that $U_c$ should be contained in the domain of definition of $\mc E$ and $\mc E$ should be $\mf X$-clc on $U_c$. In our case, the coefficients $Ri_* \mc G$ are defined on all of $X$,  and since $i=i_n\cdots i_1$, \cite[Lemma V.3.9.b]{Bo} implies that $Ri_*\mc G$ will be $\mf X$-clc on any stratification of $X$ such that $X^{n-1}=\Sigma$. \qedsymbol

In this approach, dependence of the intersection cohomology modules on the choice of $\Sigma$ occurs right at the level of coefficients. Thus with $\Sigma$ built into the ground floor, we see that the support and cosupport conditions in $A_{\bar p}$ can be relatively limited compared to those in $AX2'_{\bar p,\mf X, \mc G}$. On the other hand, the coefficients of the theory are forced to be non-clc, though they are determined from our original locally constant system $\mc G$ once we have fixed $\Sigma$. 

This conclusion also allows us to formulate the analogue of $AX2''_{\bar p,\Sigma, \mc G}$, which we will denote $AX3''_{\bar p,\Sigma, \mc G}$. These axioms, which characterize the Deligne sheaf but are independent of  the precise topological stratification subject to $\Sigma$, read as follows:
\begin{description}
\item[$(3''a)$] $\mc S^*$ is bounded, $\mc S^j=0$ for $j<0$,  $\mc S^*$ is $\mf X$-clc for some topological stratification of $X$ subject to $\Sigma$, and $\mc S^*_c=(Ri_*\mc G)|_{U_c}$ with respect to this stratification. 

\item[$(3''b)$] $\dim \{\text{supp} \mc H^j(\mc S^*)\}\leq n-\bar p^{-1}(j)$ for all $j>c-2$.
\item[$(3''c)$] for all $j<n$, $\dim \{x\in X\mid H^j(f_x^!\mc S^*)\neq 0\}\leq n-\bar q^{-1}(n-j)$.
\end{description}

It is clear that if $\mc S^*$ satisfies $AX3_{\bar p,\hat{ \mf X}, \mc G}$ for some stratification $\hat{ \mf X}$ subject to $\Sigma$ then it also satisfies $AX3''_{\bar p,\Sigma, \mc G}$. Conversely, if $\mc S^*$ satisfies $AX3''_{\bar p,\Sigma, \mc G}$ and $\hat {\mf X}$ is any topological stratification of $X$, then by Theorem \ref{T: main}, $\mc S^*$ is quasi-isomorphic to the Deligne sheaf constructed with respect to $\hat {\mf X}$. So $\mc S^*$ also satisfies $AX3_{\bar p,\hat{ \mf X}, \mc G}$.

\bibliographystyle{amsplain}
\bibliography{bib}

\end{document}